\documentclass{article}

\usepackage{a4wide}
\usepackage{amsmath}
\usepackage{amssymb}
\usepackage{amsthm}
\usepackage{txfonts}
\usepackage{mathtools}
\usepackage{tikz}
\usepackage{mdwlist}

\newtheorem{remark}{Remark}
\newtheorem{theorem}{Theorem}
\newtheorem{ass}{Assumption}

\newcommand{\pt}{\partial}
\newcommand{\laplace}{\Delta}
\newcommand{\eps}{\varepsilon}
\newcommand{\RR}{\mathbb{R}}
\renewcommand{\phi}{\varphi}
\renewcommand{\rho}{\varrho}
\renewcommand{\theta}{\vartheta}
\newcommand{\scal}[2]{\left(#1,#2\right)}
\newcommand{\dual}[2]{\left\langle#1,#2\right\rangle}
\newcommand{\bigscal}[2]{\bigl(#1,#2\bigr)}
\newcommand{\bigdual}[2]{\bigl\langle#1,#2\bigr\rangle}
\newcommand{\m}[1]{\mathcal{#1}}

\newcommand{\E}{\mathrm{e}}
\newcommand{\D}{\mathrm{  d}}

\newcommand{\abs}[1]{\left|#1\right|}
\newcommand{\norm}[1]{\left\|#1\right\|}
\newcommand{\bignorm}[1]{\big\|#1\big\|}

\title{Maximum-norm a posteriori error bounds for an extrapolated Euler/finite element
    discretisation of parabolic equations}

\author{Torsten Lin\ss\thanks{Fakult\"at f\"ur Mathematik und Informatik,
        FernUniversit\"at in Hagen,
        Universit\"atsstra{\ss}e 11,
        58095 Hagen,
        Germany,
        \texttt{torsten.linss@fernuni-hagen.de}}
   \and Goran Radojev\thanks{Department of Mathematics and Computer Science, Faculty of Sciences,
        University of Novi Sad, Trg Dositeja Obradovi\'ca~4, 21000 Novi Sad,
        Serbia,
        \texttt{goran.radojev@dmi.uns.ac.rs}.
        GR acknowledges financial support from DAAD and by FernUniversität in Hagen
        through visiting grants.}}

\begin{document}

\maketitle

\begin{abstract}
  A class of linear parabolic equations are considered.
  We give a posteriori error estimates in the maximum norm for
  a method that comprises extrapolation applied to the backward
  Euler method in time and finite element discretisations in space.
  We use the idea of elliptic reconstructions and certain bounds
  for the Green’s function of the parabolic operator.
  

\end{abstract}

\section{Introduction}

Residual-type a posteriori error estimates in the maximum norm for parabolic equations
have been given in a number of publication~\cite{MR2519598, MR2629992, MR1335652,
MR3032709, MR3056758}.

Given a second-order linear elliptic operator $\m{L}$ in a spatial domain
$\Omega\subset\RR^n$ with Lipschitz boundary, we consider the linear
parabolic equation:
\begin{subequations}\label{problem}
\begin{alignat}{2}
 \m{K}u \coloneqq \pt_t u  + \m{L} u & = f\,,
            &\quad& \text{in} \quad Q\coloneqq\Omega \times (0,T],\\
 \intertext{subject to the initial condition}
   u(x,0) & = u^0(x)\,, && \text{for} \quad x\in\bar{\Omega}, \\
 \intertext{and a homogeneous Dirichlet boundary condition}
   u(x,t) & = 0\,, && \text{for} \quad (x,t)\in \pt\Omega \times [0,T].
\end{alignat}
\end{subequations}
Precise assumptions on the data will be given later.

We consider extrapolation applied to the first-order backward Euler discretisation
in time and FEM in space applied to problem~\eqref{problem}, and obtain computable
a posteriori error estimates in the maximum norm.
The analysis follows the framework of~\cite{MR3056758}.
We also draw ideas from~\cite{MR2519598, MR2034895} and employ elliptic
reconstructions in the analysis.

The paper is organised as follows.
In Section~\ref{sect:weak-discr} we specify our assumptions on the data of
problem~\eqref{problem}, recapitulate certain aspects of the existence theory
for~\eqref{problem} and introduce our discretisation by the extrapolated Euler
method and finite elements.
In Section~\ref{sect:ana} we conduct an a posteriori error analysis of the
discretisation.
We formulate our assumptions on the existence of error estimators for the
elliptic problems, \S\ref{ssect:apost-ell}, and of certain bounds for the Green's
function of the parabolic problem, \S\ref{ssect:green}.
In \S\ref{ssect:reco} the consept of elliptic reconstructions is introduced,
while the main result, Theorem~\ref{theo:full_extra} is derived in \S\ref{ssect:apost-para}.
Finally, numerical results are presented in Section~\ref{sect:numer} to illustrate our
theoretical findings.

\paragraph{Notation.} Throughout, we denote by $\norm{\cdot}_{q,\Omega}$ the standard
norm in $L_q(\Omega)$, \mbox{$q\in[0,\infty]$}.

\section{Weak formulation and discretisation}
\label{sect:weak-discr}

We shall study~\eqref{problem} in its standard variational form,
cf.~\cite[\S5.1.1]{MR2362757}.
The appropriate Gelfand triple consists of the spaces
\begin{gather*}
  V = H_0^1(\Omega), \quad H= L_2(\Omega) \quad\text{and}\quad
  V^* = H^{-1}(\Omega)\,.
\end{gather*}
Moreover, by \mbox{$a(\cdot,\cdot) \colon V \times V \to \RR$}
we denote the bilinearform associated with the elliptic operator $\m{L}$,
while $\dual{\cdot}{\cdot}\colon V^*\times V \to \RR$ is the duality
pairing and $\scal{\cdot}{\cdot}\colon H\times H\to\RR$ is the scalar
product in~$H$.

The solution $u$ of~\eqref{problem} may be considered as a mapping
\mbox{$[0,T] \to V\colon t \mapsto u(t)$}, and we will denote its
(temporal) derivative by $u'$ (and $\pt_t u$).
Let
\begin{gather*}
  W_2^1(0,T;V,H) \coloneqq \left\{ v \in L_2(0,T;V) \colon
         v'\in L_2(0,T;V^*)\right\}\,.
\end{gather*}

Our variational formulation of~\eqref{problem} reads:
Given \mbox{$u^0\in H$} and \mbox{$F\in L_2(0,T;V^*)$},
find \mbox{$u\in W_2^2(0,T;V,H)$} such that
\begin{subequations}\label{weak}
\begin{alignat}{2}
   \frac{\D}{\D t} \bigscal{u(t)}{\chi} + a\bigscal{u(t)}{\chi}
            & = \bigdual{F(t)}{\chi} \quad \forall \chi\in V,
      \ \ t\in(0,T], \\
  \intertext{and}
    u(0)=u^0.
\end{alignat}
\end{subequations}
This problem possesses a unique solution.

In the sequell we shall assume that the source term $F$ has more
regularity and can be represented as \mbox{$\dual{F(t)}{\chi}=\scal{f}{\chi}$},
\mbox{$\forall v\in V$}, with a function \mbox{$f\in(0,T;H)$}

Since we are interested in maximum-norm error estimates we have to make further
assumptions on the data to ensure that the solution can be evaluated pointwise.
To this end, we assume that the intial and boundary data satify the zero-th
order compatibility condition, i.e. $u^0 = 0$ on $\pt\Omega$, and that $u^0$
is Hölder continuous in $\bar\Omega$.
Under standard assumptions on $f$ and $\m{L}$, problem~\eqref{problem} possesses
a unique solution that is continuous on $\bar{Q}$; see~\cite[\S5, Theorem 6.4]{0174.15403}.

Now we turn to discretising~\eqref{weak}.
To this end, let the mesh in time be given by
\begin{gather*}
  \omega_t \colon 0=t_0<t_1<\ldots<t_M=T,
  \intertext{with mesh intervals}
     I_j\coloneqq(t_{j-1},t_j) \quad
  \text{and step sizes} \quad
  \tau_j\coloneqq t_j-t_{j-1},\ j=1,2,\dots,M.
\end{gather*}
For any function \mbox{$v\colon \Omega\times[0,T]\to\RR$} that is continuous
in time on \mbox{$[0,T]$} we set
\begin{gather*}
  v^0 \coloneqq v(\cdot,0) \quad\text{and}\quad
  v^{j-\mu} \coloneqq v(\cdot,t_j-\mu \tau_j)\, \ \ j=1,2,\dots,M, \ \
    \mu\in[0,1]
\end{gather*}
Let $V_h$ be a finite dimentional (FE-)subspace of $V$ and
let $a_h\scal{\cdot}{\cdot}$ and $\scal{\cdot}{\cdot}_h$
be approximations of the bilinear form $a\scal{\cdot}{\cdot}$
and of the scalar product $\scal{\cdot}{\cdot}$ in $H$.
These may involve quadrature, for example.

Let $u^0_h\in V_h$ be an approximation of the initial condition $u^0$.
Then our discretisation of the initial-boundary-value problem~\eqref{weak}
is based on an extrapolation of the implicit Euler method and FEM in
space and reads as follows:

\begin{subequations}\label{euler-FE}
  \noindent
  \textbf{One-step Euler:}
    Set $v_h^0=u_h^0$ and find $v^j_h \in V_h$, $j=1,\dots,M$, such that
    \begin{align}\label{euler-FE-one}
       \scal{\frac{v_h^j-v_h^{j-1}}{\tau_j}}{\chi}_h + a_h \scal{v^j_h}{\chi}
          & = \scal{f^j}{\chi}_h \quad \forall \ \chi \in V_h.
  \intertext{\textbf{Two-step Euler:}
    Set $w_h^0=u_h^0$ and find $w_h^{j-1/2},w^j_h \in V_h$, $j=1,\dots,M$,
    such that}
       \label{euler-FE-two-1}
       \scal{\frac{w_h^{j-1/2} - w_h^{j-1}}{\tau_j/2}}{\chi}_h
          + a_h \scal{w_h^{j-1/2}}{\chi}
          & = \scal{f^{j-1/2}}{\chi}_h \quad \forall \ \chi \in V_h, \\
       \label{euler-FE-two-2}
       \scal{\frac{w_h^j - w_h^{j-1/2}}{\tau_j/2}}{\chi}_h
          + a_h \left( w_h^j,\chi\right)
          & = \scal{f^j}{\chi}_h \quad \forall \ \chi \in V_h.
  \intertext{\textbf{Extrapolation:} Set}
       u_h^j & \coloneqq 2 w_h^j - v_h^j, \quad j=1,\dots,M.
    \end{align}
\end{subequations}

Finally, set
\begin{gather*}
  \delta_t v^j \coloneqq \frac{v^j-v^{j-1}}{\tau_j}\,, \ \ j=1,2,\dots,M.
\end{gather*}


\section{Error analysis}
\label{sect:ana}

Our analysis of the discretisation~\eqref{euler-FE} uses three main ingredients:
\begin{itemize*}
  \item a posteriori error bounds for the elliptic problem $\m{L}y=g$,
    see \S\ref{ssect:apost-ell},
  \item bounds for the \textsc{Green}'s function associated with the
    parabolic operator $\m{K}$, see \S\ref{ssect:green} and
  \item the idea of elliptic reconstructions introduced by
    Makridakis and Nochetto~\cite{MR2034895}, see \S\ref{ssect:reco}.
\end{itemize*}
After these concepts have been reviewed, we derive an a posteriori error
bound for the extrapolated Euler method in~\S\ref{ssect:apost-para}.

\subsection{A posteriori error estimation for the elliptic problem}
\label{ssect:apost-ell}

Given $g\in H$, consider the elliptic boundary-value problem of
finding $y\in V$ such that
\begin{gather}\label{prob:ell}
   a\scal{y}{\chi} = \scal{g}{\chi}\,, \ \ \forall \chi\in V,
\end{gather}
and its discretisation of finding $y_h\in V_h$ such that
\begin{gather}\label{FEM:ell}
   a_h\scal{y_h}{\chi} = \scal{g}{\chi}_h\,, \ \ \forall \chi\in V_h.
\end{gather}

\begin{ass}\label{ass:ee}
  There exists an a posteriori error estimator $\eta$ for the
  FEM~\eqref{FEM:ell} applied to the elliptic problem~\eqref{prob:ell} with
  \begin{gather*}
    \norm{y_h-y}_{\infty,\Omega} \le \eta\bigl(y_h, g\bigr).
  \end{gather*}
\end{ass}

A few error estimators of this type are available in the literature.
We mention some of them.
\begin{itemize*}
  \item Nochetto et al.~\cite{05068766} study the semilinear problem
        $-\laplace u + g(\cdot,u)=0$ in up to three space dimensions.
        They give a posteriori error bounds for arbitrary order FEM
        on quasiuniform triangulations.
  \item Demlow \& Kopteva~\cite{zbMATH06618549} too consider arbitrary order FEM
        on quasiuniform triangulations, but for the singularly perturbed
        equation $-\eps^2\laplace u + g(\cdot,u)=0$.
        A posteriory error estimates are established that are robust in the
        perturbation parameter.
        Furthermore, in~\cite{Kop14} for the same problem $P_1$-FEM on
        \emph{anisotropic} meshes are investigated.
  \item In~\cite{MR2334045,MR3232628} arbitrary order FEM for
        the linear problem \mbox{$-\eps^2 u'' + ru = g$} in \mbox{$(0,1)$},
        \mbox{$u(0)=u(1)=0$} are considered.
        In contrast to the afore mentioned contributions all constants
        appearing in the error estimator are given explicitely.
\end{itemize*}

\subsection{Green's functions}
\label{ssect:green}

Let the \textsc{Green}'s function associated with $\m{K}$ and an arbitrary
point \mbox{$x\in\Omega$} be denoted by $\m{G}$,
Then for all \mbox{$\phi\in W^1_2\left(0,T;V,H\right)$}
\begin{gather}\label{green-rep}
  \phi(x,t) = \bigscal{\phi(0)}{\m{G}(t)}
              + \int_0^t \bigdual{\bigl(\m{K}\phi\bigr)(s)}{\m{G}(t-s)} \D s.
\end{gather}

The \textsc{Green}'s function
\mbox{$\m{G} \colon \bar{\Omega}\times[0,T]\to\RR$}, \mbox{$t \in (0,T]$}.
solves for fixed $x$
\begin{gather*}
  \pt_t\m{G}+\m{L}^*\m{G} = 0, \ \ \text{in} \ \Omega\times\RR^+, \ \
  \m{G}\bigr|_{\pt\Omega} = 0, \ \ \m{G}(0) = \delta_x=\delta(\cdot-x)\,.
\end{gather*}

\begin{ass}\label{ass:green}
  There exist non-negative constants
  $\kappa_0$, $\kappa_1$, $\kappa_1'$ and $\gamma$ such that
  \begin{gather}\label{source:ass}
    \norm{\m{G}(t)}_{1,\Omega} \le \kappa_0\,\E^{-\gamma t} \eqqcolon \phi_0(t),
       \quad
    \norm{\pt_t \m{G}(t)}_{1,\Omega}
       \le \left(\frac{\kappa_1}{t} +\kappa_1'\right)
                  \,\E^{-\gamma t} \eqqcolon \phi_1(t),
  \end{gather}
  for all $x\in\bar\Omega$, $t\in[0,T]$.
\end{ass}

In \S\ref{sect:numer} we will present numerical results for an example
test problem that satisfies these assumptions.
A more detailed discussion of problem classes for which such results are
available is given in~\cite[\S2]{MR3720388}, see also Appendix A in~\cite{MR3056758}.

\subsection{Elliptic reconstruction}
\label{ssect:reco}
Given an approximation $\phi_h^{j-\mu}\in V_h$ of $u(t_{j-\mu})$, we define
$\psi_\phi^{j-\mu}\in V_h$ by
\begin{gather}\label{reconstr-psi}
  \scal{\psi_\phi^{j-\mu}}{\chi}_h = a_h\scal{\phi_h^{j-\mu}}{\chi}
              - \scal{f^{j-\mu}}{\chi}_h \quad \forall \ \chi\in V_h\,,
      \ \ j= 0,\dots,M.
\end{gather}
This can be written as an ,,elliptic'' problem:
\begin{gather}\label{elliptic:discr}
  a_h\scal{\phi_h^{j-\mu}}{\chi} = \scal{f^{j-\mu} + \psi_\phi^{j-\mu}}{\chi}_h \quad
      \forall \ \chi\in V_h, \ \ j=0,\dots,M.
\end{gather}
Next, define $R_\phi^{j-\mu}\in H_0^1(\Omega)$ by
\begin{gather}\label{reconstr}
  a\scal{R_\phi^{j-\mu}}{\chi} = \scal{f^{j-\mu}+\psi_\phi^{j-\mu}}{\chi}
    \quad \forall \ \chi\in H_0^1(\Omega)\,, \ \ 0=1,\dots,M,
\end{gather}
or for short: $\m{L}R_\phi^{j-\mu} = f^{j-\mu}+\psi_\phi^{j-\mu}$.
The function $R_\phi^{j-\mu}$ is referred to as the elliptic reconstruction of
$\phi_h^{j-\mu}$, \cite{MR2034895}.
Later we shall employ reconstructions $R_u$, $R_v$ and $R_w$ of the approximations
$u_h$, $v_h$ and $w_h$ computed by~\eqref{euler-FE}.

Now, $\phi_h^{j-\mu}$ can be regarded as the finite-element approximation of $R_\phi^{j-\mu}$
obtained by~\eqref{elliptic:discr}, and the error can be bounded using the elliptic estimator
from~\S\ref{ssect:apost-ell}:
\begin{gather}\label{reconstr-esti}
  \norm{\phi_h^{j-\mu} - R_\phi^{j-\mu}}_\infty
    \le \eta_\mathrm{ell}^{j-\mu} \coloneqq \eta\left(\phi_h^{j-\mu}, f^{j-\mu}+\psi_\phi^{j-\mu}\right)\,,
      \ \ j=0,\dots,M.
\end{gather}
Because of linearity, we have
\begin{gather}\label{reconstr-delta-esti}
  \norm{\delta_t\left(\phi_h - R\right)^j}_\infty
    \le \eta_\mathrm{ell,\delta}^j
    \coloneqq \eta\left(\delta_t \phi_h^j, \delta_t\left(f+\psi_\phi\right)^j\right)\,,
      \ \ j=1,\dots,M.
\end{gather}

\subsection{A posteriori error estimation for the parabolic problem}
\label{ssect:apost-para}

We are now in a position to derive our a posteriori error bound
for~\eqref{euler-FE}.
We like to use the \textsc{Green}'s function representation~\eqref{green-rep}
with $\phi$ replaced by the error \mbox{$u-u_h$}.
First we have to extend the \mbox{$u_h^j$}, \mbox{$j=0,1,\dots,M$}, to a
function defined on all of $[0,T]$.
We use piecewise linear interpolation:
For any function $\phi$ defined on $\omega_t$, $t_j\mapsto \phi^j$,
we define
\begin{gather}\label{hat_notation}
  \hat{\phi}(\cdot,t) \coloneqq \frac{t_{j}-t}{\tau_j}\, \phi^{j-1}
  + \frac{t-t_{j-1}}{\tau_j}\,\phi^{j}
  \quad\text{for}\ \ t\in[t_{j-1},t_j],\quad j=1,\dots, M.
\end{gather}

Eq.~\eqref{green-rep} yields for the error at final time $T$ and for
any $x\in\Omega$:
\begin{gather}\label{rep-error}
  (u-u_h^M)(x) =
  (u-\hat{u}_h)(x,T) = \bigscal{u^0-u_h^0}{\m{G}(t)}
              + \int_0^t \bigdual{\bigl(\m{K}(u-\hat{u}_h\bigr)(s)}{\m{G}(t-s)} \D s.
\end{gather}

Next, we derive a representation of the residuum of $\hat{u}_h$ in the
differential equation.
Consider the reconstruction $R^j$ of $v_h^j$.
By~\eqref{reconstr-psi}
\begin{gather*}
  \scal{\psi_v^j}{\chi}_h = a_h\scal{v_h^j}{\chi} - \scal{f^j}{\chi}_h\,,
     \quad \forall v\in V_h, \ \ j=0,\dots,M.
\end{gather*}
Comparing with~\eqref{euler-FE-one}, we see that $\psi_v^j = \delta_tv_h^j$,
$j=1,\dots,M$.
Therefore,
\begin{gather*}
  \m{L} R_v^j = f^j -\delta_t v_h^j\,, \ \ j=1,\dots,M.
\end{gather*}
Similarly, by~\eqref{euler-FE-two-1} and~\eqref{euler-FE-two-2}
\begin{gather*}
  \m{L} R_w^{j-1/2} = f^{j-1/2} - 2 \frac{w_h^{j-1/2} - w_h^{j-1}}{\tau_j}
  \quad \text{and} \quad
  \m{L} R_w^j = f^j - 2 \frac{w_h^j - w_h^{j-1/2}}{\tau_j}\,,\ \ j=1,\dots,M.
\end{gather*}
The last three equations imply
\begin{gather*}
  \pt_t \hat{u}_h(s) = \delta_t u^j_h
        = 2 \delta_t w^j_h - \delta_t v^j_h
        = f^{j-1/2} - \m{L} \left(R_w^j + R_w^{j-1/2}-R_v^j\right)\,,
   \ \ s \in I_{j}, \ j=1,\dots,M.
\end{gather*}
For the residuum we get
\begin{align*}
  \left(\m{K} \left(u- \hat{u}_h\right)\right)(s)
     & = f(s) - \pt_t \hat{u}_{h} (s) - \left(\m{L}\hat{u}_h\right) (s) \\
     & = f(s) - f^{j-1/2} + \m{L} \left(R_w^j + R_w^{j-1/2} - R_v^j\right)
           - \left(\m{L} \bigl(\hat{u}_h-\hat{R}_u\bigr)\right) (s)
           - \left(\m{L}\hat{R}_u\right)(s)\,.
\end{align*}
For the last term on the R.H.S.,%
\footnote{Note, that
  $
    R^j_u = 2 R^j_w - R^j_v\,, \ \ j=0,1,\dots,M,
    \quad \text{and}\quad
    \hat{R}_u = 2 \hat{R}_w - \hat{R}_v \ \ \text{on} \ [0,T],
  $
  properties that will be used frequently.
}
\begin{align*}
  \left(\m{L} \hat{R}_u\right) (s)
     & = \frac{\m{L} R_u^j + \m{L} R_u^{j-1}}{2}
           + \bigl(s - t_{j-1/2} \bigr) \delta_{t} \m{L} R_u^j \\
     & = \m{L} R_w^j + \m{L} R_w^{j-1} -
       \frac{\m{L} R_v^j + \m{L} R_v^{j-1}}{2}
           + \bigl(s - t_{j-1/2} \bigr) \delta_{t} \left(\psi_u^j+f^j\right),
\end{align*}
because $\m{L}R_u^j = \delta_t\psi_u^j +f^j$.
Therefore,
\begin{align*}
  & \left(\m{K} \left(u- \hat{u}_h\right)\right)(s) \\
  & \qquad = f(s) - f^{j-1/2}
        + \m{L} \underbrace{\left(R_w^{j-1/2} - R_w^{j-1} - \frac{R_v^j - R_v^{j-1}}{2}
                \right)}_{\eqqcolon R_*^j}
        - \left(\m{L}\bigl(\hat{u}_{h}-\hat{R}_u\bigr)\right)(s)
        - \bigl(s-t_{j-1/2}\bigr) \delta_{t} \left(\psi_u^j+f^j\right).
\end{align*}
The definitions of the elliptic reconstructions $R_w$ and $R_v$ yield
\begin{gather*}
  \m{L}R_*^j = \psi_w^{j-1/2} - \psi_w^{j-1} - \frac{\psi_v^j-\psi_v^{j-1}}{2}
                   - \frac{f^j-2f^{j-1/2}+f^{j-1}}{2}
       \eqqcolon \psi_*^j - f_*^j\,, \ \ j=1,\dots,M,
\end{gather*}
respectively,
\begin{gather*}
  a\left(R_*^j,\chi\right)
     = \scal{\psi_*^j - f_*^j}{\chi} \quad \forall \ \chi\in H_0^1(\Omega)\,.
\end{gather*}
Set
\begin{gather*}
  z_*^j \coloneqq w_h^{j-1/2} - w_h^{j-1} - \frac{v_h^j - v_h^{j-1}}{2}\quad
\end{gather*}
and note, that
\begin{gather*}
  a_h\left(z_*^j,\chi\right)
     = \scal{\psi_*^j - f_*^j}{\chi}_h \quad \forall \ \chi\in V_h\,, \\
\end{gather*}
Thus, the function $z_*^j$ can be interpreted as a FE approximation of $R_*^j$,
and we have the bound
\begin{gather}\label{est-ell-z}
  \norm{R_*^j - z_*^j}_{\infty,\Omega} \le \eta\left(z_*^j, \psi_*^j-f_*^j\right)\,,
    \ \ j=1,\dots,M.
\end{gather}


Setting,
\begin{gather*}
  F(s)\coloneqq f(s)-f^{j-1/2}, \ \ \text{for} \ t\in (t_{j-1},t_j), \ \ j=1,\dots,M,
\end{gather*}
we have the following representation of the residuum:
\begin{gather*}
  \left(\m{K} \left(u- \hat{u}_h\right)\right)(s)
    = \bigl(F-\hat{F}\bigr)(s)
        - \left(\m{L}\bigl(\hat{u}_{h}-\hat{R}_u\bigr)\right)(s)
        - \bigl(s-t_{j-1/2}\bigr) \delta_{t} \psi_u^j
        + \psi_*^j - f_*^j\,, \ \ s\in I_j\,.
\end{gather*}
This is substituted into \eqref{rep-error} to obtain
\begin{gather}\label{error-rep-extra}
 \begin{split}
  & u(x,T) - u_h^M(x) \\
  & \qquad = \scal{u^0-u_h^0}{\m{G}(T)}
      + \int_0^T \scal{(F-\hat{F})(s)}{\m{G}(T-s)}\D s
      + \int_0^T \dual{\mathcal{L} (\hat{R}_u-\hat{u}_{h})  (s)}{\m{G}(T-s)} \D s \\
    & \qquad \qquad
              - \sum_{j=1}^{M}\int_{I_j} \bigl(t_{j-1/2}-s\bigr)
                \scal{\delta_{t} \psi_u^j}{\m{G}(T-s)}\D s
              + \sum_{j=1}^{M} \int_{I_j} \scal{\psi_*^j-f_*^j}{\m{G}(T-s)} \D s \,.
 \end{split}
\end{gather}

\begin{theorem}\label{theo:full_extra}
  Let  $u_h^M$ be the approximation of $u(T)$ given by~\eqref{euler-FE}.
  Then, for any $K\in\{0,\dots,M-1\}$, one has
  \begin{gather*}
    \norm{u(T)-u_h^M}_{\infty,\Omega}
      \le \eta_\mathrm{eE}^{M,K}
      \coloneqq \kappa_0 \sigma_0 \eta_{\mathrm{init}}
        + \sum_{j=1}^{M} \sigma_j
                  \left(\kappa_0 \eta_F^j+\eta_\mathrm{ell}^{M,K}
         + \chi_j \eta_{\delta\psi_u}^j+\eta_{z_h}^j\right),
  \end{gather*}
  where
  $F(s)\coloneqq f(s)-f^{j-1/2}$, for $t\in (t_{j-1},t_j)$,
  \begin{gather*}
    \eta_{\mathrm{init}}\coloneqq \norm{u^0-u_h^0}_{\infty,\Omega}, \quad
    \eta_F^j \coloneqq \int_{I_j} \norm{(F-\hat{F})(s)}_{\infty,\Omega} \D s, \quad
    \eta_{\delta\psi_u}^j \coloneqq \norm{\delta_t \psi_u^j}_{\infty,\Omega}\,. \\
    \eta_{z_h}^j \coloneqq
       \min  \Biggl\{ \kappa_0 \tau_j \bignorm{\psi_*^j-f_*^j}_{\infty,\Omega} , 
                     \mu_j \left(\bignorm{z_*^j}_{\infty,\Omega} + \eta\left(z_*^j,\phi_*^j-f_*^j\right)\right)
             \Biggr\}. \\
    \eta_\mathrm{ell}^{M,K}
        \coloneqq \kappa_0\left(\eta_\mathrm{ell}^M + \sigma_K \eta_\mathrm{ell}^K
         + \sum_{j=K+1}^M \sigma_j \tau_j \eta_{\mathrm{ell},\delta}^j \right)
         + \sum_{j=1}^K 
            \sigma_j \mu_j \max\left\{\eta_\mathrm{ell}^j, \eta_\mathrm{ell}^{j-1}\right\}\,, \\
    \sigma_j \coloneqq\E^{-\gamma(T-t_j)}, \ \
    \mu_j\coloneqq \int_{I_j} \left(\frac{\kappa_1}{T-s} + \kappa_1'\right) \D s\,, \ \
    \chi_j \coloneqq \min\left\{\frac{\kappa_0 \tau_j^2}{4},
                             \int_{I_j} \frac{(t_j-s)(s-t_{j-1})}2
                                \left(\frac{\kappa_1}{T-s} + \kappa_1'\right) \D s
               \right\}.
  \end{gather*}
  The elliptic estimators $\eta_\mathrm{ell}^j$ and
  $\eta_{\mathrm{ell},\delta}^j$ have been defined in
  \eqref{reconstr-esti} and~\eqref{reconstr-delta-esti}.
\end{theorem}
\begin{proof}
  We have to bound the right-hand side of~\eqref{error-rep-extra}
  and consider the various terms separately.
  \paragraph{(i)}
  The H\"older inequality and~\eqref{source:ass} give
  \begin{gather}\label{est-init}
    \abs{\scal{u^0-u_h^0}{\m{G}(T)}}
       \le \kappa_0 \E^{-\gamma T} \eta_{\mathrm{init}}
    \intertext{and}\label{est-F}
    \abs{\int_0^T \scal{(F-\hat{F})(s)}{\m{G}(T-s)}\D s}
       \le \kappa_0 \sum_{j=1}^M \E^{-\gamma(T-t_j)} \eta_f^j\,.
  \end{gather}

  \paragraph{(ii)}
  For the third term on the right-hand side of~\eqref{error-rep-extra},
  we have
  \begin{gather*}
    \int_0^T \scal{\m{L}\bigl(\hat{R}-\hat{u}_h\bigr)(s)}{\m{G}(T-s)}\D s
      = \int_0^T \scal{\pt_t \m{G}_t(T-s)}{\bigl(\hat{R}-\hat{u}_h\bigr)(s)}\D s,
  \end{gather*}
  because $\m{L}^*\m{G} = \pt_t \m{G}$.
  For any $K\in\{0,\dots,M-1\}$, integration by parts on $(t_K,T)$ , gives
  \begin{align*}
    & \int_0^T \scal{\pt_t\m{G}(T-s)}{\bigl(\hat{R}-\hat{u}_h\bigr)(s)}\D s \\
    & \qquad
       = - \scal{\m{G}(0)}{\bigl(R-u_h\bigr)^M}
         + \scal{\m{G}(T-t_K)}{\bigl(R-u_h\bigr)^K}
         + \sum_{j=K+1}^M \int_{I_j}
           \scal{\m{G}(T-s)}{\delta_t\bigl(R-u_h\bigr)^j} \D s \\
    & \qquad\qquad
         + \sum_{j=1}^K \int_{I_j} \scal{\pt_t\m{G}(T-s)}{\bigl(\hat{R}-\hat{u}_h\bigr)(s)}\D s
  \end{align*}
  We apply H\"older's inequality, \eqref{source:ass},
  \eqref{reconstr-esti} and~\eqref{reconstr-delta-esti}
  to obtain
  \begin{gather}\label{est-L(R-u)}
    \begin{split}
    & \abs{\int_0^T \dual{\m{L}\bigl(\hat{R}-\hat{u}_h\bigr)(s)}{\m{G}(T-s)}\D s} \\
    & \qquad
       \le \kappa_0\left(\eta_\mathrm{ell}^M + \E^{-\gamma(T-t_K)} \eta_\mathrm{ell}^K
         + \sum_{j=K+1}^M \E^{-\gamma(T-t_j)} \tau_j \eta_{\mathrm{ell},\delta}^j \right)
         + \sum_{j=1}^K
            \int_{I_j} \phi_1(T-s) \D s \
                   \max\left\{\eta_\mathrm{ell}^j, \eta_\mathrm{ell}^{j-1}\right\}\,.
    \end{split}
  \end{gather}

  \paragraph{(iii)}
  The fourth term in~\eqref{error-rep-extra} is bounded as follows.
  \begin{gather*}
    \abs{\int_{I_j} \bigl(t_{j-1/2}-s\bigr)
                  \scal{\m{G}(T-s)}{\delta_t\psi_u^j} \D s}
        \le \frac{\kappa_0 \tau_j^2}{4} \E^{-\gamma(T-t_j)}
            \norm{\delta_t\psi_u^j}_{\infty,\Omega}
  \end{gather*}
  An alternative bound is obtained using integration by parts.
  Let $\omega(s)=\frac{1}{2}\bigl(t_j-s\bigr)\bigl(s-t_{j-1}\bigr)$, $s\in[t_{j-1},t_j]$,
  and note that $\frac{\D}{\D s}\omega(s) = t_{j-1/2}-s$.
  Then, we have
  \begin{gather*}
    \int_{I_j} \bigl(t_{j-1/2}-s\bigr) \scal{\m{G}(T-s)}{\delta_t\psi_u^j} \D s
      = \int_{I_j} \omega(s) \scal{\pt_t\m{G}(T-s)}{\delta_t\psi_u^j}\D s\,,
  \end{gather*}
  and estimate as follows
  \begin{gather*}
    \abs{\int_{I_j} \bigl(t_{j-1/2}-s\bigr)
                  \scal{\m{G}(T-s)}{\delta_t\psi_u^j} \D s}
        \le \int_{I_j} \omega(s)
                     \left(\frac{\kappa_1}{T-s}+\kappa_1'\right) \D s \
                     \E^{-\gamma(T-t_j)} \norm{\delta_t\psi_u^j}_{\infty,\Omega}\,,
  \end{gather*}
  by~\eqref{source:ass}, $p=1$.
  Combining these two bounds, we get
  \begin{gather}\label{est-dpsi}
  \abs{\int_{I_j} \bigl(t_{j-1/2}-s\bigr)
                  \scal{\m{G}(T-s)}{\delta_t \psi_u^j}\D s}
                     \leq \E^{-\gamma(T-t_j)} \chi_j
           \bignorm{\delta_t \psi_u^j}_{\infty,\Omega} \E^{-\gamma(T-t_j)}\,.
  \end{gather}

  \paragraph{(iv)}
  For the last term in~\eqref{error-rep-extra} we proceed as follows, again
  using H\"older's inequality and~\eqref{source:ass}.
  \begin{gather*}
    \abs{\int_{I_j} \scal{\m{G}(T-s)}{\psi_*^j-f_*^j}\D s}
        \leq  \kappa_0 \tau_j \E^{-\gamma(T-t_j)}\bignorm{\psi_*^j-f_*^j}_{\infty,\Omega}\,.
  \end{gather*}
  Furthermore,
  \begin{gather*}
    \scal{\m{G}(T-s)}{\psi_*^j-f_*^j} = 
    \dual{\m{L}R_*^j}{\m{G}(T-s)} = 
    \scal{\pt_t\m{G}(T-s)}{R_*^j-z_*^j} + \scal{\pt_t\m{G}(T-s)}{z_*^j}\,,
  \end{gather*}
  which provides a second bound:
  \begin{gather*}
    \abs{\int_{I_j} \scal{\m{G}(T-s)}{\psi_*^j - f_*^j}\D s}
      \le \int_{I_j} \phi_1(T-s) \D s \
          \left\{ \bignorm{R_*^j-z_*^j}_{\infty,\Omega}
                 +\bignorm{z_*^j}_{\infty,\Omega}
          \right\}.
  \end{gather*}
  Combining both bounds, we get
  \begin{gather}\label{est-z}
    \begin{split}
      &\abs{\int_{I_j}  \scal{\m{G}(T-s)}{\psi_*^j - f_*^j}\D s} \\
      & \qquad \quad
      \le \E^{-\gamma(T-t_j)} \min  \left\{ \kappa_0 \tau_j  \bignorm{\psi_*^j-f_*^j}_{\infty,\Omega}, 
           \int_{I_j} \left(\frac{\kappa_1}{T-s}+\kappa_1'\right) \D s
                 \left(\bignorm{z_*^j}_{\infty,\Omega} + \eta\left(z_*^j,\psi_*^j-f_*^j\right)\right)  \right\}.
    \end{split}
  \end{gather}

  Finally, applying~\eqref{est-init}--\eqref{est-z} to~\eqref{error-rep-extra}
  completes the proof.
\end{proof}

\begin{remark}\label{rem:est}
  \emph{(i)}
  In general, the supremum norm involved in $\eta_\mathrm{init}$ can not be
  determined exactly, but needs to be approximated.
  For example, one can use a mesh that is finer than the finite-element mesh.

  \emph{(ii)}
  The integral in $\eta_F^j$ needs to be approximated.
  One possibility is  Simpson's rule, which is of higher order
  and gives
  \begin{gather*}
    \int_{I_j} \norm{\bigl(F-\hat{F}\bigr)(s)}_{\infty,\Omega} \D s
       \approx \frac{\tau_j}{6} \norm{f^j - 2 f^{j-1/2} + f^{j-1}}_{\infty,\Omega}
       \,.
  \end{gather*}
  Here too, the supremum norm needs to be approximated.
\end{remark}


\section{A numerical example}
\label{sect:numer}

Consider the following reaction-diffusion equation
\begin{subequations}\label{testproblem}
\begin{alignat}{2}
 \pt_t u  - u_{xx} + (5x+6) u & = \E^{-4t} - \cos(x+t)^4\,,
            &\quad& \text{in} \quad (-1,1) \times (0,1],\\
 \intertext{subject to the initial condition}
   u(x,0) & = \sin \frac{\pi(1+x)}{2}\,, && \text{for} \quad x\in[-1,1], \\
 \intertext{and the Dirichlet boundary condition}
   u(x,t) & = 0\,, && \text{for} \quad (x,t)\in \{-1,1\} \times [0,1].
\end{alignat}
\end{subequations}
The Green's function for this problem satisfies
\begin{gather*}
  \norm{\m{G}(t)}_{1,\Omega} \le \E^{-t/2},
    \quad
  \norm{\pt_t\m{G}(t)}_{1,\Omega}
            \le \frac{3}{2^{3/2}} \frac{\E^{-t/2}}{t} \,, \ \
    \text{see~\cite{MR3720388}}.
\end{gather*}
The exact solution to this problem is unknown.
To compute a reference solution, we use a spectral method in space combined
with the dG(2) method in time which is of order $5$. This gives an approximation
that is accurate up to machine precision.

Our spatial discretisation uses the version of $P_1$-FEM analysed
in~\cite{MR2334045} and the a posteriori estimator derived therein.
The method is of order $2$, and we couple spetial and temporal mesh
sizes by \mbox{$h=\tau$}.

Table~\ref{Table-full_ext} displays the results of our test computations.
The first column contains the number of mesh intervals $M$ (with \mbox{$h=\tau=1/M$}),
followed by the errors $e_M$ at final time, the experimental order of
convergence $p_M$, the error estimator $\eta_\mathrm{eE}^{M,M-1}$ and finally
the efficiency $\chi_M$:
\begin{gather*}
  e_M\coloneqq \bignorm{u(T)-U^M}_{\infty,\Omega}, \quad
  p_M\coloneqq \frac{\ln (e_{M/2}/e_M)}{\ln 2} \quad \text{and} \quad
  \chi_M \coloneqq \frac{\eta_\mathrm{eE}^{M,M-1}}{e_M}.
\end{gather*}
\begin{table}[!h]
\begin{center}
\begin{tabular}{ccccc}
$M$ & $e_M$ & $p_M$ & $\eta_\mathrm{eE}^{M,M-1}$ & $\chi_M$ \\ \hline
    $2^4$ & 3.872e-04 & 1.90 & 4.038e-01 & 1/1043 \\
    $2^5$ & 1.039e-04 & 1.94 & 1.050e-01 & 1/1011  \\
    $2^6$ & 2.703e-05 & 1.97 & 2.647e-02 & 1/979  \\
    $2^7$ & 6.908e-06 & 1.99 & 6.646e-03 & 1/962 \\
    $2^8$ & 1.742e-06 & 2.00 & 1.667e-03 & 1/957  \\
  $2^{9}$ & 4.369e-07 & 2.00 & 4.175e-04 & 1/956 \\
  $2^{10}$ & 1.092e-07 & 2.00 & 1.045e-04 & 1/957 \\
  $2^{11}$ & 2.730e-08 & 2.00 & 2.617e-05 & 1/958 \\
  $2^{12}$ & 6.824e-09 & 2.00 & 6.549e-06 & 1/960  \\
  $2^{13}$ & 1.706e-09 & 2.00 & 1.639e-06 & 1/961 \\
  $2^{14}$ & 4.301e-10 & 1.99 & 4.102e-07 & 1/954 \\
 \hline
\end{tabular}
\caption{Error, estimator and efficiency, test problem~\eqref{testproblem}\label{Table-full_ext}}
\end{center}
\end{table}
The numbers confirm our finding in Theorem~\ref{theo:full_extra}.
The errors are overestimated by a factor of about $1000$.

Table~\ref{Table-full_ext2} displays the various components of the error estimator
from Theorem~\ref{theo:full_extra}.
The dominant term is $\eta_\mathrm{ell}^{M,M-1}$, which contains the contributions
from the elliptic error estimator.

\begin{table}[!h]
\begin{center}
\centerline{
\begin{tabular}{cccccccccc}
$M$ &  $\eta_{\mathrm{init}}$ & $\eta_F$ & $\eta_\mathrm{ell}^{M,M-1}$
    & $\eta_{\delta\psi_u}$ & $\eta_{z_h}$ \\ \hline
    $2^4$  & 5.696e-04  & 1.418e-02  & 3.628e-01  & 4.379e-03 & 2.186e-02 \\
    $2^5$  & 1.425e-04  & 3.522e-03  & 9.623e-02  & 1.244e-03 & 3.855e-03 \\
    $2^6$  & 3.564e-05  & 8.706e-04  & 2.445e-02  & 3.428e-04 & 7.684e-04 \\
   $2^7$   & 8.910e-06  & 2.162e-04  & 6.151e-03  & 9.279e-05 & 1.777e-04 \\
   $2^8$   & 2.228e-06  & 5.387e-05  & 1.542e-03  & 2.488e-05 & 4.430e-05 \\
   $2^9$   & 5.569e-07  & 1.344e-05  & 3.858e-04  & 6.629e-06 & 1.106e-05 \\
  $2^{10}$ & 1.392e-07  & 3.358e-06  & 9.651e-05  & 1.758e-06 & 2.774e-06  \\
  $2^{11}$ & 3.481e-08  & 8.391e-07  & 2.413e-05  & 4.646e-07  & 6.949e-07  \\
  $2^{12}$ &  8.702e-09  & 2.097e-07  & 6.034e-06  & 1.224e-07  & 1.742e-07  \\
  $2^{13}$ & 2.175e-09  & 5.243e-08  & 1.509e-06  & 3.216e-08  & 4.366e-08  \\
 $2^{14}$  & 5.439e-10  & 1.311e-08  & 3.772e-07  & 8.431e-09  & 1.096e-08  \\
 \hline
\end{tabular}
}
\caption{Composition of the error estimator, test problem~\eqref{testproblem}\label{Table-full_ext2}}
\end{center}
\end{table}

\bibliographystyle{plain}

\def\cprime{$'$}

\end{document}